\documentclass{article} 

\usepackage{amsfonts,amsmath}
\usepackage{amssymb,latexsym}
\usepackage[dvips]{epsfig}
\usepackage{amscd}  

\title{\bf A category for the adjoint representation}

 \author{Ruth Stella Huerfano\footnote{Departamento de Matem\'{a}ticas, 
Universidad National de Colombia, Santaf\'{e} de Bogot\'{a}, Colombia,\newline 
huerfano@matematicas.unal.edu.co}   
\hspace{0.05in} 
and Mikhail Khovanov\footnote{Department of Mathematics, University 
of California, Davis, mikhail@math.ucdavis.edu}} 

\date{}

\newtheorem{prop}{Proposition}

\newtheorem{lemma}{Lemma}
\newtheorem{corollary}{Corollary}

\newcommand{\oplusop}[1]{{\mathop{\oplus}\limits_{#1}}}

\begin{document}

\maketitle

\baselineskip 12pt

%
%

\def\C{\mathbb C}
\def\R{\mathbb R}
\def\N{\mathbb N}
\def\Z{\mathbb Z}
\def\Q{\mathbb Q}
\def\F{\mathbb F}
\def\P{\mathbb P}
\def\SS{\mathbb S}
\def\Zq{\Z[q,q^{-1}]}
\def\l{\lbrace}
\def\r{\rbrace}
\def\o{\otimes}
\def\D{\Delta}
\def\O{\mathcal{O}}
\def\lra{\longrightarrow}
\def\slt{\mathfrak{sl}_2}
\def\U{\mathbb U}
\def\Ug{U_q(\mathfrak g)}
\def\dU{\dot{\mathbb U}} 
\def\Qq{\Q(q)}

\def\Br{\mathrm{Br}}
\def\a{\mathbf a} 
\def\l{\mathbf l} 
\def\H{\mathbb H} 
\def\M{\mc{M}}
\def\mC{\mc{C}}
\def\A{\mathcal{A}}
\def\Gr{\mbox{Gr}}
\def\G0{\Gamma^0}
\def\Bred{B^{\mathrm{red}}}
\def\Hom{\mbox{Hom}}
\def\det{\mbox{det}}
\def\tE{\widetilde{E}}
\def\tF{\widetilde{F}}
\def\wt{\widetilde}

\def\Ea{\mc{E}_{\alpha}}
\def\Fa{\mc{F}_{\alpha}}
\def\Ca{\mc{C}_{\alpha}}

\def\binom#1#2{\left( \begin{array}{c} #1 \\ #2 \end{array}\right)}
\def\qbinom#1#2{\left[ \begin{array}{c} #1 \\ #2 \end{array}\right]}

\def\doublemaprights#1#2#3#4{\raise3pt\hbox{$\mathop{\,\,\hbox to
     #1pt{\rightarrowfill}\kern-30pt\lower3.95pt\hbox to
     #2pt{\leftarrowfill}\,\,}\limits_{#3}^{#4}$}}

\def\mc{\mathcal} 
\def\mf{\mathfrak}
\def\cC{{\mathcal{C}}}
\def\cR{{\mathcal{R}}}

\def\CV{{\C\mbox{-Vect}}}

\def\sl{\mathfrak{sl}}
\def\TP1{T^{\ast}\P^1} 

\def\yesnocases#1#2#3#4{\left\{
\begin{array}{ll} #1 & #2 \\ #3 & #4 
\end{array} \right. }

\def\drawing#1{\begin{center} \epsfig{file=pictures/#1} \end{center}}

\def\hsm{\hspace{0.05in}}


\section{Introduction} 
\label{intro} 

The adjoint representation of a simple Lie algebra $\mf{g}$ admits a 
deformation into an irreducible representation $R$ of the quantum group 
$U_q(\mf{g}).$ In this paper for a simply-laced $\mf{g}$ we realize $R$ 
as the Grothendieck group of a particular abelian category $\cC.$ There 
are exact functors from $\cC$ to $\cC$ which on the Grothendieck group 
act as the quantum group generators $E_{\alpha}, F_{\alpha},$ where  
$\alpha$ varies over simple roots. Various relations in the 
quantum group between products 
of $E_{\alpha}$ and $F_{\alpha}$ become functor isomorphisms. 

The adjoint representation $R$ has a weight space decomposition as the 
direct sum of
$1$-dimensional vector spaces, one for each  root of $\mf{g},$ and 
the Cartan subalgebra.  Mirroring this, 
we define $\cC$ as the direct sum of copies of the category of graded 
vector spaces and the category of graded modules over the 
algebra $A(\Gamma),$ naturally associated to the Dynkin 
diagram $\Gamma$ of $\mf{g}.$ Change each edge 
of $\Gamma$ into a pair of oriented edges, form the path algebra of this 
oriented graph, and quotient it out by the ideal generated 
by certain linear combinations of length $2$ paths. $A(\Gamma)$ 
is the resulting quotient algebra, and we name it the \emph{zigzag algebra} 
of $\Gamma.$ The Grothendieck group of the category 
of $A(\Gamma)$-modules is naturally identified with the weight lattice 
in the Cartan subalgebra of $\mf{g}.$ 

We introduce functors $\mc{E}_{\alpha}$ and $\mc{F}_{\alpha}$ lifting the 
generators $E_{\alpha}$ and $F_{\alpha}$ of $U_q(\mf{g})$
and check that defining relations in the quantum group become 
isomorphisms of functors. We proceed to explore various properties 
of our categorification of the quantum group action on $R.$ 
Among them is the 
adjointness of functors $\mc{E}_{\alpha}$ and 
$\mc{F}_{\alpha},$ existence of several dualities in $\cC$ and 
a braid group action in 
the derived category of $A(\Gamma)$-modules. 

We expect that not just the adjoint but 
any finite-dimensional 
irreducible representation $L$ of the quantum group $U_q(\mf{g}),$ for 
a simple simply-laced Lie algebra $\mf{g},$ admits a canonical realization 
as the Grothendieck group of an abelian category $\cC(L).$ 
In this realization 
the Kashiwara-Lusztig basis in $L$ should become the basis 
of indecomposable projective objects, the quantum group should  
act by exact functors and there should be a braid group action in the derived 
category of $\cC(L).$ In short, all structures of the category $\cC$ that 
we describe in this paper should also be present in categories $\cC(L).$  
Categories $\cC(L)$ will be very close relatives of categories of 
coherent sheaves on Nakajima quiver varieties [Na] and categories of 
modules over cyclotomic Hecke algebras [A]. 
The  
work of Ariki [A], among other things, contains a categorification 
of all irreducible finite-dimensional 
representations of $\mf{sl}_n.$ His
categories are made of blocks of the categories of modules 
over cyclotomic Hecke algebras for generic $q.$ 
Ariki's goals, which include a proof and generalizations of the 
Lascoux-Leclerc-Thibon 
conjecture [LLT], are quite different from ours. In particular, it 
has not been checked whether various fine structures of the category $\cC,$
described in Section~\ref{cat} of our paper and expected to hold 
in categories $\cC(L),$ are present in Ariki's categories. 

This work is intended to provide a  simple model example of a "perfect"
categorification, with all structures visible in the 
representation $L$ lifted to its categorification $\cC(L).$ 
Another model example, a categorification of irreducible 
$U_q(\mf{sl}_2)$ representations, will be treated in [Kh]. 

Our second goal is to draw the reader's attention to the zigzag algebra 
$A(\Gamma)$ of a graph $\Gamma.$ 
Zigzag algebras have a variety of nice features, which we 
discuss in Sections~\ref{zzalgebras} and \ref{McKay}:  

(i) $A(\Gamma)$ is a trivial extension algebra and 
has a nondegenerate symmetric trace form;  

(ii) if $\Gamma$ is a finite Dynkin diagram, 
then $A(\Gamma)$ has finite type and 
there is a bijection between indecomposable representations 
of $A(\Gamma)$ and roots of $\mf{g};$ 

(iii) if $\Gamma$ is bipartite, the quadratic dual of the zigzag algebra 
is the preprojective algebra of $\Gamma,$ for a sink-source orientation of 
$\Gamma;$  

(iv) Any multiplicity one Brauer tree algebra is derived equivalent 
to a zigzag algebra; 

(v) if $\Gamma$ is a bipartite affine Dynkin diagram, $A(\Gamma)$ is Morita 
equivalent to the cross-product algebra $\Lambda(\C^2,G),$ where 
$G$ is the finite subgroup of $SU(2)$ associated to $\Gamma$ via the 
McKay correspondence, and $\Lambda(\C^2,G)$ 
is the cross-product of the group algebra of $G$ and 
the exterior algebra on $2$ generators. 

Section~\ref{McKay} expands on (v) to explain the role played by 
the zigzag algebras in the McKay correspondence. 
This section can be viewed as a comment on a recent work of 
Kapranov and Vasserot [KV], where categories of coherent sheaves 
on resolutions of simple surface singularities are related to 
categories of 
modules over the cross-product of $\C[x,y]$ and the group algebra 
of $G,$ the latter cross-product 
being Koszul dual to $\Lambda(\C^2, G).$ 

We conclude the paper with 
Section~\ref{rep-theory}, where we compile a  
surprisingly long and diverse list of other appearances of zigzag 
algebras in the representation theory and geometry.

\vspace{0.1in}

\emph{Acknowledgements} We are very much indebted to the referee for the 
Journal of Algebra for pointing out a number of errors in the first 
version of this paper. We are grateful to Paul Seidel 
for pointing out Proposition~\ref{triv-ext} and for interesting 
discussions. This paper naturally branched out of the joint work [KS] of 
Paul Seidel and the second author. 

During our work on the paper 
M.K. was partially supported by NSF grants DMS 9729992 and DMS 9627351.

\section{The adjoint representation of a simply-laced quantum group} 
\label{the-adjoint} 

\subsection{Quantum groups} 
\label{qgroups} 

 Let $\mf{g}$ be a complex simple simply-laced Lie algebra, $\Phi$ 
the root system of $\mf{g}$ and  $\Pi$ a set of simple roots. The Weyl 
group  $W$ of $\mf{g}$ acts on the real vector space $\R\Phi$ and there is a 
unique $W$-invariant bilinear form on $\R\Phi$ such that 
  $(\alpha, \alpha)=2$ for any root $\alpha\in \Phi.$ 


Let $\Q(q)$ be the field of polynomial functions with rational coefficients 
in an indeterminate $q.$ 

The quantum group $U=U_q(\mf{g})$ is a $\Q(q)$-algebra with 
generators $E_{\alpha}, F_{\alpha}, K_{\alpha},$  $ K^{-1}_{\alpha}$ 
for $\alpha\in \Pi$   and relations  
\begin{equation}
\label{relaciones}    
\begin{split} 
&   K_{\alpha} K^{-1}_{\alpha} = 1 = K^{-1}_{\alpha} K_{\alpha}, \\
&   K_{\alpha} K_{\beta} = K_{\beta} K_{\alpha}, \\
&   K_{\alpha} E_{\beta} = q^{(\alpha, \beta)} E_{\beta}K_{\alpha}, \\
&   K_{\alpha} F_{\beta} = q^{-(\alpha, \beta)} F_{\beta} K_{\alpha}, \\
&   E_{\alpha} F_{\beta} - F_{\beta}E_{\alpha} = \delta_{\alpha \beta}
  \frac{K_{\alpha} - K_{\alpha}^{-1}}{q-q^{-1}}, \\
&   E_{\alpha} E_{\beta} = E_{\beta}E_{\alpha} \mbox{ for }
  (\alpha,\beta) = 0, \\  
&   F_{\alpha} F_{\beta} = F_{\beta}F_{\alpha} \mbox{ for }
  (\alpha,\beta) = 0, \\  
&   E_{\alpha}^2 E_{\beta} - (q+q^{-1}) E_{\alpha} E_{\beta} E_{\alpha} 
  + E_{\beta} E_{\alpha}^2 = 0 \mbox{ for } 
  (\alpha, \beta) = -1, \\
&   F_{\alpha}^2 F_{\beta} - (q+q^{-1}) F_{\alpha} F_{\beta} F_{\alpha} 
  + F_{\beta} F_{\alpha}^2 = 0 \mbox{ for } 
  (\alpha, \beta) = -1. 
\end{split}
\end{equation}  
Let $\overline{\phantom{a}}$ be the $\Q$-linear involution of $\Q(q)$ 
which changes $q$ into $q^{-1}.$ 

$U$ has an antiautomorphism $\tau: U \to U^{\mathrm{op}}$ described by 
\begin{equation} 
\begin{split}   
& \tau(E_{\alpha})=q F_{\alpha}K_{\alpha}^{-1}, \hsm 
\tau(F_{\alpha})= q E_{\alpha}K_{\alpha}, \hsm 
\tau(K_{\alpha}) =K_{\alpha}^{-1},   \\
& \tau(fx) = \overline{f}\tau(x), \hsm \mbox{ for }f\in \Q(q)
\mbox{ and } x\in U,   \\
& \tau(xy) = \tau(y)\tau(x), \hsm \mbox{ for }x,y\in U. 
 \label{welcome-tau}
\end{split} 
\end{equation} 
Let $\psi$ be the $\Q$-algebra involution of $U$ defined by 
\begin{equation}
\label{vote-for-psi}
\begin{split}   
& \psi(E_{\alpha})= E_{\alpha}, \hsm \psi(F_{\alpha}) = F_{\alpha}, 
\hsm \psi(K_{\alpha}) = K_{\alpha}^{-1}, \\
& \psi (fx) = \overline{f} x \hsm \mbox{ for } \hsm f\in \Q(q) 
\hsm \mbox{ and } \hsm x\in U. 
\end{split} 
\end{equation} 
Let $\omega$ be the $\Q(q)$-linear involution of $U$ given by 
\begin{equation} 
\label{for-omega} 
\omega(E_{\alpha}) = F_{\alpha}, \hsm \omega(F_{\alpha}) = E_{\alpha}, 
\hsm \omega(K_{\alpha})= K_{\alpha}^{-1}
\end{equation}
 
These three automorphisms and antiautomorphisms $\tau, \psi, \omega$ 
satisfy the following relations
\begin{equation} 
\label{commute} 
 \psi \omega = \omega\psi, \hsm  \tau \omega = \omega \tau, 
  \hsm   \tau \psi \tau = \psi.  
\end{equation}

Given a $\Qq$-vector space $V,$ a $\Q$-bilinear form $V\times V\to \Q(q)$ 
is called semilinear if it is 
$\Q(q)$-antilinear in the first variable and $\Q(q)$-linear in the 
second, i.e.
\begin{equation}
\label{semilin}  
<fx,y>= \overline{f}<x,y>, \hsm <x,fy>=f<x,y>\hsm\mbox{ for } f\in \Qq
\hsm \mbox{ and } x,y\in U. 
\end{equation}

  \subsection{The adjoint representation} 
  \label{adj-rep} 

  The adjoint representation of the quantum group $U_q(\mf{g})$ is 
  the irreducible representation with the highest weight equal to the 
  maximal root. Denote this representation by $R.$ It has a basis 
  $ \{  x_{\mu}, h_{\alpha} \}$   for $\mu \in \Phi, \alpha\in \Pi$  
  with the following action of the quantum group:  
  \begin{equation}
 \label{Kalpha} 
  K_{\alpha} x_{\mu} = q^{(\alpha, \mu)} x_{\mu}, 
  \hspace{0.2in} 
  K_{\alpha} h_{\beta} = h_{\beta};  
  \end{equation} 
$E_{\alpha}, F_{\alpha}$ act by 
\begin{equation} 
   \label{action-EF}
   \begin{array}{lll} 
   E_{\alpha} x_{\mu} = 0,  & 
   F_{\alpha} x_{\mu} = 0   & 
   \mbox{if }(\mu,\alpha) = 0; \\ 
   E_{\alpha} x_{\mu} = 0,  & 
   F_{\alpha} x_{\mu} = x_{\mu-\alpha}   & 
   \mbox{if }(\mu,\alpha) = 1; \\ 
   E_{\alpha} x_{\mu} = x_{\mu+\alpha},  & 
   F_{\alpha} x_{\mu} = 0   & 
   \mbox{if }(\mu,\alpha) = -1; \\
   E_{\alpha} x_{\alpha} = 0,  & 
   F_{\alpha} x_{\alpha} = h_{\alpha},  & \\
   E_{\alpha} x_{-\alpha} = h_{\alpha}, & 
   F_{\alpha} x_{-\alpha} = 0; & 
   \end{array}   
   \end{equation} 
 and, for $\alpha, \beta \in \Pi,$ 
\begin{equation} 
\label{action-EF-2}    
\begin{array}{lll} 
   E_{\alpha} h_{\beta} & = &  \left\{ 
   \begin{array}{ll}
    (q+q^{-1}) x_{\alpha}  & \mbox{ if }\alpha=\beta \\
   x_{\alpha} & \mbox{ if }(\beta,\alpha) 
   = -1 \\
   0   & \mbox{ otherwise}
   \end{array} 
  \right.   \\
   F_{\alpha} h_{\beta} & = &  
   \left\{ 
   \begin{array}{ll} (q+q^{-1}) x_{-\alpha}  & \mbox{ if }\alpha=\beta \\
                         x_{-\alpha} & \mbox{ if }(\beta,\alpha) 
   = -1 \\
   0   & \mbox{ otherwise}
   \end{array} 
  \right. 
  \end{array}
\end{equation} 
   
We will denote by $R_{\mu}$ the weight $\mu$ subspace of $R.$ For 
$\mu\in \Phi,$ the subspace $R_{\mu}$ is one-dimensional, while the 
dimension of $R_0$ is equal to the rank of $\mf{g}.$ 

   On $R$ there is a $\tau$-invariant semilinear form $<,>,$  
\begin{equation} 
\label{tau-invariant} 
 <xa,b> = <a,\tau(x)b> \mbox{ for any }x\in U\mbox{ and } a,b\in R, 
\end{equation} 
described in our basis by 
\begin{equation} 
\label{product-semi} 
\begin{array}{ll}  
   <x_{\mu}, x_{\mu}>= 1, &  \mu \in \Psi \\
   <h_{\alpha}, h_{\alpha}>= 1+ q^2, &   \alpha\in \Pi \\
   <h_{\alpha}, h_{\beta}>= q  &  \mbox{if }(\alpha, \beta)= -1, 
   \hspace{0.1in} \alpha, \beta \in \Pi \\
   <h_{\alpha}, h_{\beta}>= 0  &  \mbox{if }(\alpha, \beta)= 0 , 
   \hspace{0.1in} \alpha, \beta \in \Pi 
\end{array} 
\end{equation}  
Distinct weight components of $R$ are orthogonal with respect to this 
form. 
The basis $\{ x_{\mu}, h_{\alpha}\}$ is called the canonical basis of $R.$
It is a special case of the Lusztig-Kashiwara basis [L1],[L2],[Ka] 
(also see [J] for an introduction) in 
irreducible $U_q(\mf{g})$ representations. 
  
Let $l_{\alpha}\in R_0$ be defined by $< h_{\beta}, l_{\alpha}> = 
\delta_{\alpha, \beta}.$ The basis $\{ x_{\mu}, l_{\alpha}\}$ is dual 
to the canonical basis of $R$ with respect to the semilinear form $<,>.$ 
We call this basis the dual canonical basis of $R.$ 

Denote by $I$ the $\Z[q,q^{-1}]$-submodule of $R$ generated by elements 
of the dual canonical basis and by $I'$ the $\Z[q,q^{-1}]$-submodule of 
$R$ generated by canonical basis vectors. Note that 
$I'$ is a $\Z[q,q^{-1}]$-submodule
of $I.$

Let $\psi_{R}$ be the $\Q$-linear involution $R\to R$ given by 
\begin{equation} 
\label{psiR} 
\psi_R(x_{\mu}) = x_{\mu}, \hsm\hsm  
\psi_R(h_{\alpha}) = h_{\alpha}, \hsm \hsm 
\psi_R(fv) = \overline{f}\psi_R(v)\hsm \mbox{ for } f\in \Q(q), v\in U. 
\end{equation} 
It is clear that $\psi_R(ax) = \psi(a) \psi_R(x)$ for $a\in U$ and $x\in R$ 
and that 
\begin{equation}
\label{psiR-inv}  
<\psi_R x, \psi_R y > = < y, x> \hspace{0.1in} \mbox{ for } \hspace{0.1in} 
x,y\in R. 
\end{equation} 
Let $\omega_R$ 
be the $\Q(q)$-linear involution $R\to R$ which takes $x_{\mu}$ to 
$x_{-\mu}$ and $h_{\alpha}$ to $h_{\alpha}.$ We have 
$\omega_R(ax) = \omega(a) \omega_R(x)$ for $a\in U$ and $x\in R,$ and 
\begin{equation} 
\label{u-invariance} 
<\omega_R x, \omega_R y> 
= <x,y>\hspace{0.1in} \mbox{ for }\hspace{0.1in} x,y\in R.
\end{equation} 
 Note that involutions $\psi_R $ and $\omega_R$ preserve $I\subset R.$

\section{Graphs and algebras} 
\label{gr-and-alg} 

Let $\Gamma$ be a connected graph without loops and multiple edges. 
Associated to $\Gamma$ there is the double graph, $D\Gamma,$ which 
has the same vertices as $\Gamma$ and twice as many edges as $\Gamma.$ 
Namely, each edge $f$ of $\Gamma$ is substituted by two oriented edges, 
which connect the same vertices as $f$ and have opposite 
orientations. This construction is best illustrated by an example:

\begin{center} \epsfig{file=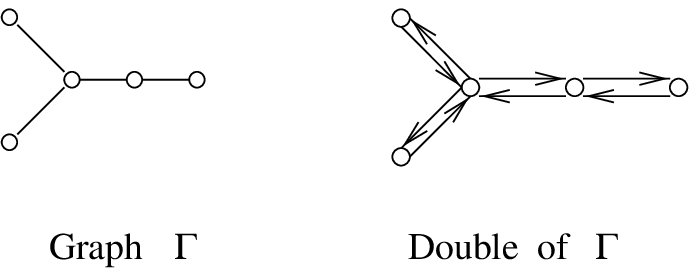} 
\end{center} 

Take the path algebra of $D\Gamma.$ It is an algebra (over $\C$) 
spanned by all oriented paths in $D\Gamma$ with the multiplication 
given by concatenating paths. In particular, minimal idempotents 
correspond one-to-one to length $0$ paths, i.e., to vertices of 
$\Gamma.$ Since we assume that $\Gamma$ has no multiple edges, we 
can describe a path by a list of vertices it travels through, thus, a 
path that starts at a point $a$, goes to $b$ and then to $c$ will 
be denoted $(a|b|c).$  If $\Gamma$ has more than 2 vertices, 
denote by $A(\Gamma)$ the quotient algebra of this path 
algebra by the ideal generated by the following elements

(i) Paths $(a|b|c)$ for each triple of vertices $a,b,c$ of $\Gamma$ such that  
$a,b$ are connected, $b,c$ are connected and $a\not=c,$ 

(ii) Element $(a|b|a)-(a|c|a)$ whenever $a$ is connected to both $b$ and 
$c.$ 

If $\Gamma$ consists of the single vertex
only, define $A(\Gamma)$ as the algebra generated by $1$ and $X$ with 
$X^2=0.$ If $\Gamma$ consists of two points joined by a single edge, 
define $A(\Gamma)$ as the quotient of the path algebra of $D\Gamma$  
by the two-sided ideal spanned by all paths of length greater than $2.$ 
We will call $A(\Gamma)$ the {\it zigzag algebra} of $\Gamma.$ 

If $\Gamma$ has more than one vertex, 
$A(\Gamma)$ has a natural grading with paths of length $k$ in
 degree $k.$ If $\Gamma$ have only one vertex, we introduce a 
grading on  $A(\Gamma)= \C[X]/(X^2)$ by placing $X$ in degree $2.$ 

Denote by $v(\Gamma)$ the set of vertices of $\Gamma$ and by $e(\Gamma)$ 
the set of its edges. 
For any $\Gamma,$ the algebra $A(\Gamma)$  has nonzero components only 
in degrees $0,1$ and $2,$ and the dimensions of these graded components 
are $v(\Gamma), 2 e(\Gamma) $ and $v(\Gamma),$ respectively. 

\begin{prop} 
\label{sym-algebra} 
$A(\Gamma)$ is a graded symmetric algebra. 
\end{prop} 

\emph{Proof} A finite-dimensional $\C$-algebra $A$ is called symmetric 
if it possesses a nondegenerate symmetric trace map $A\to \C.$ 
The trace map $\mbox{tr}: A(\Gamma)\to \C$ is 
defined by sending each path of length $2$ to $1$ and all other paths to $0.$ 
Clearly, this map is symmetric, $\mbox{tr}(xy) = \mbox{tr}(yx)$ 
for all $x,y\in A(\Gamma),$ and nondegenerate. $\square$ 

In particular, $A(\Gamma)$ is self-injective, i.e. $A(\Gamma)$ is injective 
as a left and right module over itself. 

\vspace{0.1in}

Let $\CV$ be the category of finite-dimensional graded complex vector spaces. 
The morphisms in this category are grading-preserving linear maps. 
If $A$ is a finite-dimensional graded $\C$-algebra, denote by 
\mbox{$A$-Mod} the abelian category of finite-dimensional graded 
$A$-modules and grading-preserving homomorphisms. 
If $M= \oplusop{n} M_n$ 
is a graded module over a graded algebra $A,$ denote by 
$M\{ k\}$ the module $M$ with the grading shifted up by $k,$ so that 
$M\{ k\}_n = M_{n-k}.$ Denote by $\{ k\}$ the functor of shifting the 
grading up by $k.$

\vspace{0.1in}

For a vertex $a\in v(\Gamma)$ denote by $P_a$ the left projective 
module $A(\Gamma)e_a$ and by $_aP$ the right projective module $e_a A(\Gamma),$
where $e_a=(a)$ is the minimal idempotent equal to the zero length 
path which begins and ends in $a.$ The left module $P_a$ is spanned by all 
paths ending in $a$ and $_aP$ is spanned by all paths starting at $a.$ 
Any indecomposable graded projective left  $A(\Gamma)$ module is 
isomorphic, up to a shift in the grading, to $P_a$ for some vertex $a.$ 
We have 
\begin{equation} 
\label{theformula} 
 _aP \o_{A(\Gamma)} P_b \cong
 \left \{ 
  \begin{array}{ll}
   \C \oplus \C\{ 2\} & \mbox{ if }\alpha=\beta \\
   \C\{ 1\}  & \mbox{ if }(\beta,\alpha) 
   = -1 \\
   0   & \mbox{ otherwise}
   \end{array} 
  \right.  
\end{equation}

For a vertex $a$ of $\Gamma$ consider functors
\begin{equation} 
\label{for-the-lemma} 
\begin{array}{ll} 
T_a: A(\Gamma)\mbox{-Mod} \lra \CV, &  
T_a(M) = _aP\o_{A(\Gamma)} M, \\
S_a: \CV \lra A(\Gamma)\mbox{-Mod}, & 
S_a(V) =  P_a \o_{\C} V. 
\end{array} 
\end{equation} 

\begin{lemma} \label{two-adjoint} 
Functor $T_a$ is right adjoint to $S_a$ and left adjoint to $S_a\{ -2\}.$ 
\end{lemma} 

This lemma follows from Proposition~\ref{sym-algebra}. Since the trace 
map of $A(\Gamma)$ has degree $-2,$ this accounts for the appearance 
of the shift by $\{ -2\}$ in the lemma. $\square$

\section{Categorification} 
\label{cat}

\subsection{Grothendieck groups} 
\label{prelim} 

If $\mc{B}$ is an abelian category, denote by $G(\mc{B})$ the Grothendieck 
group of $\mc{B}.$ The Grothendieck group 
 is the abelian group generated by symbols $[M]$ as 
$M$ ranges over all objects of $\mc{B},$ with relations $[M_2]=[M_1]+[M_3]$ 
for each short exact sequence $0\to M_1\to M_2\to M_3 \to 0.$ 
An exact functor between abelian categories 
induces a homomorphism of their Grothendieck groups. 
Denote by $G'(\mc{B})$ the subgroup of $G(\mc{B})$ generated by symbols 
$[P]$ for all projective objects $P$ of $\mc{B}.$ We call $G'(\mc{B})$ the 
projective Grothendieck group.

If $A$ is a graded algebra,  
the Grothendieck group and the projective Grothendieck group of 
the category $A$-Mod
are $\Z[q,q^{-1}]$-modules, where multiplication by $q$ corresponds to 
shift in the grading: $[M\{ k\}]= q^k [M].$

\subsection{The category $\cC$} 
\label{the-category} 

We follow the notations from Section~\ref{the-adjoint}: 
$\mf{g}$ is a  simple simply-laced Lie algebra, $\Phi$ a root system of 
$\mf{g}$ and $\Pi$ a set of simple roots. 

Let $\cC_0$ be the category $A(\Gamma)$-Mod 
of finite dimensional graded left $A(\Gamma)$-modules, 
where $\Gamma$ is the Dynkin diagram of $\mf{g}.$ 
For each root $\mu$ of the root system $\Phi,$ denote by 
$\cC_{\mu}$ the category $\CV.$ For each $\mu,$ choose a one-dimensional 
vector space in $\cC_{\mu},$ concentrated in degree $0,$ and denote it 
by $\C_{\mu}.$ 

Let the category $\cC$ be the direct sum of $\cC_0$ and  
categories $\cC_{\mu}$ for all $\mu\in \Phi$: 
\begin{equation} 
\label{cat-C} 
\cC = \oplusop{\mu\in \Psi \cup \{ 0\}} \cC_{\mu}
\end{equation} 

The vertices of the Dynkin diagram 
 $\Gamma$ are enumerated by the set of simple roots $\Pi,$ thus 
for each simple root $\alpha$ there is an associated projective 
$A(\Gamma)$-module $P_{\alpha}$ and its simple quotient $L_{\alpha}.$ 

There is an isomorphism between the Grothendieck group of 
$\cC$ and the $\Z[q,q^{-1}]$-submodule $I$ of $R,$ given by sending 
$[\C_{\mu}]\in G(\cC_{\mu})$ 
to $x_{\mu}$ and $[L_{\alpha}]$ to $q l_{\alpha}.$ 
Thus, we identify images of simple objects of $\cC$ in the Grothendieck group 
$G(\cC)$ with the elements of the dual canonical basis of $R.$ 
Under this identification indecomposable projective modules are mapped  
to the (shifted by $q$) 
canonical basis vectors:  $[P_{\alpha}] \longmapsto q h_{\alpha}.$ 
We denote this  isomorphism $G(\cC) \cong I$ by $\iota$ and 
will use it to identify the Grothendieck group $G(\cC)$ with 
$I,$ thus, we will write $[L_{\alpha}]= q l_{\alpha}, [P_{\alpha}]= 
qh_{\alpha},$ 
etc. This isomorphism restricts to an isomorphism between the projective 
Grothendieck group of $\cC$ and the submodule $I'$ of $I.$ 

The semilinear form $<,>$ on $R$ can be interpreted as dimensions of 
homomorphism spaces. Namely, if $P\in \cC$ is projective and $M\in \cC$ is 
any module, we have 
\begin{equation}
\label{Euler-char}  
\sum_{i\in \Z} q^i\mbox{dim}(\mbox{Hom}_{\cC}(P\{ i\} , M)) = < [P], [M]>. 
\end{equation} 

\subsection{Functors} 
\label{functors} 

On the category $\cC$ we define functors $\mc{E}_{\alpha}, \mc{F}_{\alpha}$ 
as follows. If $M\in \cC_0$ then 
\begin{equation} 
\label{EF-define-1} 
\begin{array}{l}  
\Ea (M) = (_{\alpha}P\o_{A(\Gamma)} M)\o \C_{\alpha}  \\
\Fa (M) = (_{\alpha}P\o_{A(\Gamma)} M) \o \C_{-\alpha} 
\end{array}
\end{equation}  
Here and further the tensor products are over $\C$ unless indicated 
otherwise. 

If $M\in \cC_{\mu}$ and $\mu \not=0$ then 
\begin{equation} 
\label{EF-define-2} 
\begin{array}{lll} 
\Ea(M)  = 0, &  \Fa(M) =0  & \mbox{ if } \hsm (\mu,\alpha) = 0 \\
\Ea(M) = 0, & \Fa(M) = M\o \C_{\mu-\alpha} & \mbox{ if } \hsm 
(\mu,\alpha) = 1 \\
\Ea(M) = M \o \C_{\mu+\alpha}, &  \Fa(M) = 0 & \mbox{ if }\hsm 
(\mu,\alpha) 
= -1 \\
\Ea(M)=0, &  \Fa(M) = P_{\alpha} \o M\{ -1\}  &  \mbox{ if } \hsm 
\mu = \alpha \\
\Ea(M) = P_{\alpha} \o M\{ -1\}, &  \Fa(M) = 0  &  \mbox{if } \hsm 
\mu = -\alpha 
\end{array}
\end{equation}  
Note that $\Ea(\cC_{\mu}) \subset \cC_{\mu+ \alpha}$ if 
$\mu+ \alpha\in \Phi \cup \{ 0 \} $ and $\Ea(\cC_{\mu}) = 0$ otherwise. 
Similarly, $\Fa(\cC_{\mu})\subset \cC_{\mu- \alpha}$ if 
$\mu - \alpha \in \Phi\cup \{ 0 \}$ and $\Fa(\cC_{\mu})= 0$ otherwise. 

Introduce the functor $\mc{K}_{\alpha}: \cC \to \cC$ by 
\begin{equation}
\label{functor-K}  
\mc{K}_{\alpha} (M) = M\{ (\mu, \alpha)\}\hspace{0.1in} \mbox{ for } 
\hspace{0.1in} M\in \cC_{\mu}
\end{equation} 
Denote by $\mc{K}^{-1}_{\alpha}$ the inverse functor to $\mc{K}_{\alpha},$ 
thus, $\mc{K}^{-1}_{\alpha}(M) = M\{ - (\mu,\alpha)\}$ for $M\in \cC_{\mu}.$ 

Earlier we identified the Grothendieck group of $\cC$ with the 
$\Z[q,q^{-1}]$-submodule $I$ of $R.$ The functors $\Ea, \Fa,\mc{K}_{\alpha}$ 
are exact, commute with the shift functor $\{ 1\},$ and on the Grothendieck 
group of $\cC$ act as the generators $E_{\alpha}, F_{\alpha}, K_{\alpha}$ 
of $U.$ 

\subsection{Quantum group relations} 
\label{liftrelations} 

\begin{prop} There are functor isomorphisms
\begin{equation} 
\label{functor-isom} 
\begin{array}{l}  
\mc{K}_{\alpha} \mc{K}^{-1}_{\alpha} \cong \mathrm{Id} 
\cong  \mc{K}^{-1}_{\alpha} \mc{K}_{\alpha}, \\ 
\mc{K}_{\alpha} \mc{K}_{\beta} \cong \mc{K}_{\beta} \mc{K}_{\alpha}, \\
\mc{K}_{\alpha} \mc{E}_{\beta} \cong
\mc{E}_{\beta}\mc{K}_{\alpha}\{(\alpha, \beta)\}, \\ 
\mc{K}_{\alpha} \mc{F}_{\beta} \cong
\mc{F}_{\beta} \mc{K}_{\alpha}\{-(\alpha, \beta)\}, \\
\mc{E}_{\alpha} \mc{F}_{\beta} \cong  \mc{F}_{\beta}\mc{E}_{\alpha}
\hspace{0.1in} \mbox{ if } \hspace{0.1in} \alpha \not= \beta, \\
\mc{E}_{\alpha} \mc{E}_{\beta} \cong \mc{E}_{\beta}\mc{E}_{\alpha}
\hspace{0.1in}\mbox{ if }\hspace{0.1in}(\alpha,\beta) = 0, \\ 
\mc{F}_{\alpha} \mc{F}_{\beta} \cong \mc{F}_{\beta}\mc{F}_{\alpha} 
\hspace{0.1in}\mbox{ if }\hspace{0.1in}(\alpha,\beta) = 0, \\  
\mc{E}_{\alpha}^2 \mc{E}_{\beta} \oplus 
 \mc{E}_{\beta} \mc{E}_{\alpha}^2 \cong 
(\mathrm{Id}\{ 1\} \oplus \mathrm{Id}\{ -1\} ) \mc{E}_{\alpha} 
\mc{E}_{\beta} \mc{E}_{\alpha} 
\hspace{0.1in} \mbox{ if } \hspace{0.1in}  
(\alpha, \beta) = -1 \\
\mc{F}_{\alpha}^2 \mc{F}_{\beta} \oplus 
 \mc{F}_{\beta} \mc{F}_{\alpha}^2 \cong 
(\mathrm{Id}\{ 1\} \oplus \mathrm{Id}\{ -1\} ) \mc{F}_{\alpha} 
\mc{F}_{\beta} \mc{F}_{\alpha} 
\hspace{0.1in} \mbox{ if } \hspace{0.1in}  
(\alpha, \beta) = -1 
\end{array} 
\end{equation} 
\end{prop} 

\emph{Proof: } In a semisimple $\C$-linear 
category any linear relation in the endomorphism algebra 
of the Grothendieck group can be lifted into a functor isomorphism. 
Our category $\cC$ has a huge semisimple direct summand, consisting of 
 categories 
$\cC_{\mu}$ for $\mu\not= 0.$ Restricted to this subcategory, functor 
isomorphisms of the above proposition  
exist for obvious reasons. To prove functor 
isomorphisms (\ref{functor-isom}) when the source category is $\cC_{\mu}$ 
for some root $\mu$ and the target category is $\cC_0,$ it is 
enough to check, for each equation in (\ref{functor-isom}), 
that the functors on the 
left and right hand sides of it, applied 
to the simple object $\C_{\mu}$, produce isomorphic objects. Since 
the target object is always projective, and isomorphism classes of 
projectives are determined by their images on the Grothendieck group, 
the claim follows. Functor isomorphisms (\ref{functor-isom}) in the 
case when $\cC_0$ is the source category and $\cC_{\mu}$ the target 
category are proved similarly, by observing that each functor is 
isomorphic to the functor of tensoring with a graded right projective 
$A(\Gamma)$-module. $\square$ 

\vspace{0.1in} 

This takes care of all defining relations (\ref{relaciones}), 
save the following one:  
\begin{equation} 
\label{lastequation} 
E_{\alpha} F_{\alpha}- F_{\alpha} E_{\alpha} = \frac{K_{\alpha} -
K_{\alpha}^{-1}}{q - q^{-1}}.
\end{equation} 
The right hand side of (\ref{lastequation}) acts by 
$\frac{q^{(\mu,\alpha)}- q^{-(\mu, \alpha)}}{ q- q^{-1}}$ on the 
weight subspace $R_{\mu}$ of $R.$ If $i = (\mu, \alpha) $ is nonnegative, 
this quotient equals to $[i] = q^{i-1} + q^{i-3} + \dots + q^{1-i},$ and is 
a Laurent polynomial in $q$ with positive coefficients. Thus, on a  weight 
subspace $R_{\mu}$ for $(\mu, \alpha) \ge 0,$ we can 
rewrite (\ref{lastequation}) as 
\begin{equation} 
E_{\alpha} F_{\alpha} = F_{\alpha} E_{\alpha} + [(\mu,\alpha)]
\end{equation} 
 Similarly, on the weight subspace $R_{\mu}$ for 
$(\mu, \alpha) \le 0$ we can rewrite the equation (\ref{lastequation}) as 
\begin{equation} 
E_{\alpha} F_{\alpha} + [-(\mu, \alpha)]= F_{\alpha} E_{\alpha}. 
\end{equation} 
Both left and right hand sides of the two equations above 
 have only positive coefficients, and it is 
in this form that the equation (\ref{lastequation}) lifts into a 
functor isomorphism. 
To state the isomorphism, we will 
denote by  $\mbox{Id}^{[j]}$ for $j\ge 0$ the functor 
$\mbox{Id}\{ j-1\} \oplus \mbox{Id}\{ j-3\} \oplus \dots \oplus \mbox{Id} 
\{ 1-j\}$ in the category $\cC_{\mu}.$ 

\begin{prop} 
\label{ef-fe-iso} 
For $\mu \in \Phi \cup \{ 0\}$ there is an isomorphism 
of functors in the category $\cC_{\mu}$ 
\begin{eqnarray} 
& & \Ea \Fa \cong \Fa \Ea \oplus \mathrm{Id}^{[(\mu, \alpha)]} 
\hspace{0.1in} \mbox{ if } \hspace{0.1in} (\mu, \alpha)\ge 0, \\
& & \Ea \Fa \oplus \mathrm{Id}^{[-(\mu, \alpha)]} \cong 
\Fa \Ea 
\hspace{0.1in} \mbox{ if } \hspace{0.1in} (\mu, \alpha)\le 0. 
\end{eqnarray} 
\end{prop} 

Proof is left to the reader. $\square$ 

\vspace{0.1in} 

In the next two subsections we show that the three automorphisms and 
antiautomorphisms of $U$, defined in Section~\ref{qgroups}, 
can be interpreted as 
various dualities in the category $\cC.$ We then explain how the 
braid group action on the weight $0$ subspace $R_0$ of $R$ lifts to 
a braid group action on the derived category of $A(\Gamma)$-modules. 
This plentitude of interesting  structures in 
the category $\cC$ clearly points to the 
naturality and uniqueness of $\cC.$ Any other 
realization of the adjoint representation $R$ as the Grothendieck group 
of an abelian category will fail to be as rich 
as the one that we describe here. 

\subsection{Adjointness and dualities}  
\label{a-and-d} 

{\bf Adjointness and the antiautomorphism} $\tau.$ Functors $\Ea$ and $\Fa$ 
have left and right adjoints: 

\begin{prop} \label{all-adjoint} 
The functor $\Ea$ is left adjoint to $\Fa\mc{K}_{\alpha}^{-1} 
\{ 1\},$ the functor $\Fa$ is left adjoint to $\Ea\mc{K}_{\alpha} 
\{ 1\}$ and $\mc{K}_{\alpha}$ is left adjoint to $\mc{K}_{\alpha}^{-1}.$ 
\end{prop} 

\emph{Proof } 
This proposition easily  reduces to Lemma~\ref{two-adjoint}.  $\square$ 

Comparing this proposition with the formula (\ref{welcome-tau}) 
for the antiautomorphism $\tau,$ we see that $\tau$ {\it becomes the 
operation of taking the right adjoint functor.} 
Let $V^{\mathrm{ad}}$ denote the right adjoint of a functor $V,$ when 
the right adjoint exists. Suppose that functors $V_1$ and $V_2$ are composable 
and admit right adjoints. Then $(V_1 V_2)^{\mathrm{ad}}\cong V_2^{\mathrm{ad}}
V_1^{\mathrm{ad}},$ i.e. passing to adjoints interchanges the order in 
the product of functors. Correspondingly, $\tau$ is an antiautomorphism, 
$\tau(ab) = \tau(b) \tau(a).$ 

Formula (\ref{tau-invariant}) now has a nice categorical interpretation. 
Namely, earlier we found that the semilinear form $<,>$ computes the 
dimension of the spaces of morphisms from a projective to an arbitrary 
module in $\cC$ (see formula (\ref{Euler-char})). If $V$ denotes an arbitrary 
product  of functors $\Ea, \Fa, \mc{K}_{\alpha}^{\pm 1}$ and $\{\pm 1\},$
there is an isomorphism 
\begin{equation} 
\label{adjointness-isomorphism} 
\mbox{Hom}_{\cC}( V (P), M ) \cong \mbox{Hom}_{\cC}(P, V^{\mathrm{ad}} (M))
\end{equation} 
for any objects $P,M$ of $\cC.$ For what follows, take $P$ a projective 
module. Then $V(P)$ is also projective, so that, after an appropriate 
summation over all shifts by $\{ i\},$ as in  
(\ref{Euler-char}), we derive the formula (\ref{tau-invariant}) for 
the $\tau$-invariance of the semilinear form $<,>.$ 

 \vspace{0.1in} 

$\psi$ {\bf and the contravariant duality functor.} 
Let $\ast: \C\mbox{-Vect}\to \C\mbox{-Vect}$ be the contravariant 
functor in the category of graded vector spaces which 
takes a vector space $V$ to its dual $V^{\ast}= \Hom(V,\C).$ 
Note that $(\C\{i\})^{\ast}\cong \C\{ -i\}.$ 

Let $\chi$ be the antiinvolution of $A(\Gamma)$ which takes 
a path $(a_1|a_2| \dots| a_k)$ in $A(\Gamma)$ to the opposite
path $(a_k|\dots|a_2|a_1).$ 

Let $\Psi: \cC \to \cC$ be the following contravariant duality 
functor. For $M\in \mbox{Ob}(\cC_{\mu})$ where $\mu\not=0$ we set 
$\Psi M = M^{\ast}\in \cC_{\mu}.$ 
For $M\in \cC_0= A(\Gamma)\mbox{-Mod}$ 
the graded vector space $M^{\ast}$ has the structure 
of a right graded $A(\Gamma)$-module. We use the antiinvolution 
$\chi$ to make it into a left graded $A(\Gamma)$-module. Thus, 
$\Psi M = M^{\ast} \in \cC_0.$ 

The following proposition is obvious. 

\begin{prop} 
\begin{enumerate} 
\item On the Grothendieck group of $\cC$ functor $\Psi$ acts as the 
involution $\psi_R$ of $R,$ defined by the formula (\ref{psiR}).  
\item There are functor isomorphisms 
\begin{equation}
\label{psi-isomo} 
\Psi \Ea \cong \Ea \Psi, \hsm \Psi \Fa \cong \Fa \Psi, \hsm 
\Psi K_{\alpha} \cong K_{\alpha}^{-1} \Psi, \hsm 
\Psi \{ i\} \cong \{ -i\} \Psi. 
\end{equation} 
\item $\Psi$ is an involution, i.e. $\Psi^2$ is isomorphic to 
the identity functor. 
\end{enumerate} 
\end{prop} 
Note that functor isomorphisms (\ref{psi-isomo}) 
lift the defining relations (\ref{vote-for-psi}) of $\psi.$ 
Thus, $\Psi$ lifts the involution $\psi$ of $U$ as 
well as the involution $\psi_R$ of $R.$ Since $\Psi$ is a contravariant 
equivalence of $\cC,$ there is an isomorphism 
\begin{equation} 
\label{iso-psi} 
\mbox{Hom}_{\cC}(\Psi M, \Psi N) \cong \mbox{Hom}_{\cC}(N,M), 
\hspace{0.2in}\mbox{for }\hspace{0.15in} M,N\in\mbox{Ob}(\cC), 
\end{equation} 
which in  the Grothendieck group of $\mc{C}$ translates into 
the formula (\ref{psiR-inv}).

\vspace{0.1in} 

$\omega_R$ {\bf as an automorphism of $\cC.$} 
Let $\Omega$ be the following self-equivalence of the category $\cC$: 

(i) $\Omega,$ restricted to $\cC_0,$ is the identity functor, 

(ii) $\Omega,$ restricted to $\cC_{\mu}$ for $\mu\in \Phi,$ is an 
equivalence of categories $\cC_{\mu} \to \cC_{-\mu},$ coming from the 
identification of both $\cC_{\mu}$ and $\cC_{-\mu}$ with the category 
$\C$-Vect. 

The following proposition is obvious. 

\begin{prop} 
\begin{enumerate} 
\item On the Grothendieck group of $\cC$ functor $\Omega$ acts as the 
involution 
$\omega_R$ of $R$ ($\omega_R$ was defined in Section~\ref{adj-rep}).  
\item There are functor isomorphisms 
\begin{equation}
\label{this-is-for-omega} 
\Omega \Ea \cong \Fa \Omega, \hspace{0.1in} 
\Omega\Fa \cong \Ea \Omega, \hspace{0.1in} 
\Omega\mc{K}_{\alpha} \cong \mc{K}_{\alpha}^{-1} \Omega. 
\end{equation} 
\item $\Omega$ is an involution, i.e. $\Omega^2$ is isomorphic to 
the identity functor. 
\end{enumerate} 
\end{prop} 

Functor $\Omega$ corresponds to the involution $\omega$ of $U$ and 
$\omega_R$ of $R.$ In particular, 
formula (\ref{for-omega}) becomes functor isomorphisms 
(\ref{this-is-for-omega}). 
Since $\Omega$ is an equivalence, there is an isomorphism 
\begin{equation} 
\label{iso-U} 
\mbox{Hom}_{\cC}(\Omega M, \Omega N) \cong \mbox{Hom}_{\cC}(M,N), 
\hspace{0.2in} M,N\in\mbox{Ob}(\cC), 
\end{equation} 
which in  the Grothendieck group of $\mc{C}$ translates into 
the formula (\ref{u-invariance}). 

\subsection{The braid group action} 
\label{braid-act} 

For a graph $\Gamma$ denote by $\Br(\Gamma)$ the braid group  
associated to $\Gamma.$ It has generators $\sigma_a$ for each vertex 
$a$ of $\Gamma$ and relations 
\begin{eqnarray*} 
\sigma_a\sigma_b\sigma_a & = & \sigma_b\sigma_a\sigma_b 
\hspace{0.2in} \mbox{ if $a$ and $b$ are joined by an edge}, \\
\sigma_a\sigma_b & = & \sigma_b\sigma_a 
\hspace{0.2in} \mbox{ otherwise.} 
\end{eqnarray*} 

Every finite dimensional representation $V$ of $U_q(\mf{g})$ comes equipped 
with a natural  action of the braid group $\Br(\Gamma),$ 
where $\Gamma$ is the Dynkin diagram of $\mf{g}$ (see Jantzen [J], for 
instance). This action preserves the weight $0$ subspace of $V.$ When 
$V=R$ is the adjoint representation, the braid group action on the weight 
$0$ subspace $R_0$ has a particularly simple form, with $\sigma_{\alpha}$ 
acting as $qE_{\alpha} F_{\alpha} - 1.$ This action can be 
categorified. Indeed, the functor counterpart of the operator $qE_{\alpha} 
F_{\alpha}$ is $\Ea \Fa \{ 1\}.$ Restricted to the subcategory 
$\cC_0$ of $\cC$ this functor is isomorphic to the functor of tensoring 
with the $A(\Gamma)$-bimodule $P_{\alpha}\o \hsm _{\alpha}P.$ 
There is a canonical bimodule map 
\begin{equation} 
\zeta\hspace{0.1in} : \hspace{0.1in} P_{\alpha}\o \hsm _{\alpha}P 
\lra  A(\Gamma), \hspace{0.3in} 
\zeta(l_1\o l_2) = l_1 l_2 
\end{equation} 
where $l_1\in P_{\alpha}$ is a path with target $\alpha$ and $l_2\in \hsm 
_{\alpha}P$ is a path with source $\alpha.$ 

Form the bounded derived category $D^b(A(\Gamma)\mbox{-Mod})$ of the 
abelian category $A(\Gamma)$-Mod. Let 
\begin{equation} 
\Sigma_a \hspace{0.1in} : \hspace{0.1in} 
D^b(A(\Gamma)\mbox{-Mod})\lra D^b(A(\Gamma)\mbox{-Mod})
\end{equation}  
be the functor of tensoring with the complex 
\begin{equation} 
\label{braid-complex} 
0 \lra P_a\o \hsm _aP \stackrel{\zeta}{\lra} 
A(\Gamma) \lra 0
\end{equation} 
of $A(\Gamma)$-bimodules, where $a$ is a vertex of $\Gamma.$ 

\begin{prop} \label{br-action}
Functors $\Sigma_a$ are invertible  and there are 
functor isomorphisms 
\begin{eqnarray*} 
\Sigma_a\Sigma_b\Sigma_a & \cong  & \Sigma_b\Sigma_a\Sigma_b 
\hspace{0.2in} \mbox{ if $a$ and $b$ are joined by an edge}, \\
\Sigma_a\Sigma_b & \cong & \Sigma_b\Sigma_a 
\hspace{0.4in} \mbox{ otherwise.} 
\end{eqnarray*} 
\end{prop} 

Proofs in [KS] and [ST] for the case when $\Gamma$ is 
the Dynkin diagram of $\mf{sl}_n$ generalize to arbitrary graphs 
without difficulty. When the graph is a chain, this braid group action 
was also considered by R.Rouquier and A.Zimmermann~\cite{RZ}. 

 $\square$ 

Thus, the braid group $\Br(\Gamma)$ acts on the derived category of 
the category of graded $A(\Gamma)$-modules. 
When $\Gamma$ is a finite Dynkin graph, we obtain a braid group action 
on the derived category $D^b(\cC_0)$ which lifts the braid group action 
on the weight $0$ subspace of the adjoint representation. 
It is  proved in [KS] that this action is faithful if $\Gamma$ is 
the Dynkin  diagram of $\mf{sl}_n.$

Braid group and invertibility 
relations come from homotopy equivalences between 
complexes of bimodules which are tensor products of complexes 
(\ref{braid-complex}) and of similar complexes describing the inverse 
functors. In particular, functors of tensoring with 
(\ref{braid-complex}) also define a braid group action in various homotopy 
categories of complexes of $A(\Gamma)$-modules. 

Note that the 
braid group action on the $W$-invariant subspace $\oplus_{\mu \in\Psi} 
R_{\mu}$ of $R$ can be trivially lifted to the action on the derived category 
$D^b(\oplus_{\mu\in \Psi}\cC_{\mu}),$ the latter 
action given by permutations of  categories 
$\cC_{\mu}$ and shifts in the derived category. Thus, the 
braid group acts in the derived category $D^b(\cC).$ On the Grothendieck 
group this action 
descends to the standard action of the braid group in the adjoint 
representation $R$ of $U_q(\mf{g}).$ 

\vspace{0.1in}

Our braid group actions generalize to deformations of zigzag algebras 
with type (ii) (see Section~\ref{gr-and-alg}) 
defining relations $(a|b|a)=\nu_{b,c}^a (a|c|a),$ 
where $\nu_{b,c}^a\in \C^{\ast}$ and satisfy compatibility relations 
$\nu_{b,c}^a \nu_{c,b}^a=1$ and $\nu_{b,c}^a\nu_{c,d}^a
\nu_{d,b}^a=1$ whenever $a$ is connected to $b,c,$ and $d.$ 
In a zigzag algebra $\nu_{b,c}^a=1$ for all possible triples $(a,b,c).$  

We call these deformations \emph{skew-zigzag algebras} and denote by 
$A_{\nu}(\Gamma).$ They are Frobenius but not, in general, symmetric 
algebras. 
Rescaling $(a|b)$'s changes coefficients $\nu_{b,c}^a,$ and the 
moduli space of skew-zigzag algebras is naturally isomorphic to 
$H^1(\Gamma,\C^{\ast}).$ In particular, if $\Gamma$ is a tree, all of 
its skew-zigzag algebras are isomorphic. 

Derived and homotopy categories of modules over skew-zigzag algebras 
admit braid group actions that are constructed in the same way as for 
zigzag algebras.

\section{Zigzag algebras and their representations} 
\label{zzalgebras}

Let $B$ be a finite dimensional algebra over $\C.$ There is 
an obvious $B$-bimodule structure on  $B^{\ast}= \Hom_{\C}(B, \C).$  
Define the algebra $T(B)= B\oplus B^{\ast}$ with the multiplication 
$(x,f)(y,g) = (xy, fy + xg)$ for $x,y\in B$ and $f,g\in B^{\ast}.$ 
Then $T(B)$ is an associative algebra with a nondegenerate symmetric
trace map $\mbox{tr}: T(B) \to \C.$ It is called {\it the trivial extension 
algebra of } $B.$ 

Let $\Gamma$ be a graph as before and denote by $\Gamma^0$ 
any of the oriented graphs obtained by picking an orientation of each 
edge of $\Gamma.$ Let $B(\Gamma^0)$ be the path algebra of $\Gamma^0.$ 
A theorem of Gabriel says that $B(\Gamma^0)$ has finite representation 
type if and only if $\Gamma$ is a finite Dynkin diagram, in which 
case there is a natural one-to-one correspondence between indecomposable 
representations of $B(\G0)$ and positive roots of the root system associated 
to $\Gamma.$ 

Let  $\Bred(\Gamma^0)$ be the path algebra of $\Gamma^0,$ 
quotiented out by all paths of length greater than $1.$ 
Let $\Gamma^1$ be the graph $\Gamma$ with the orientation opposite to that 
of $\Gamma^0.$ Algebras $B(\Gamma^0)$ and $\Bred(\G0)$ are graded by 
the length of paths and we have the following obvious 
(here and further we refer the reader to [BGS] for definition and properties 
of Koszul algebras) 

\begin{prop} 
Algebras $\Bred(\G0)$ and  $B(\Gamma^1)$ are Koszul dual. $\square$ 
\end{prop} 

There is a natural inclusion of algebras 
$\Bred(\G0)\hookrightarrow A(\Gamma).$ 
Indeed, minimal idempotents of $\Bred(\G0)$ and of $A(\Gamma)$ are identified 
with vertices of $\Gamma$ and every oriented edge of $\G0$ is also an 
edge of the double $D\Gamma$ of $\Gamma.$ This correspondence extends 
to the abovementioned inclusion $\Bred(\G0)\hookrightarrow A(\Gamma).$ 
Note that as a vector space $A(\Gamma)$ decomposes into the direct sum of 
$\Bred(\G0)$ and the subspace spanned by edges of $D\Gamma$ which are 
complementary to the ones of $\G0$ and by length $2$ paths $(a|b|a),$
one for each vertex $a$ of $\Gamma.$ This complementary subspace, 
considered as a $\Bred(\G0)$-bimodule, is canonically isomorphic to 
$(\Bred(\G0))^{\ast},$ and, therefore, we derive 

\begin{prop}
\label{triv-ext}
 For any orientation $\Gamma^0$ of the graph $\Gamma,$ the 
algebra $A(\Gamma)$ is isomorphic to the trivial extension algebra of 
$\Bred(\G0).$ 
\end{prop}

Assume now that $\Gamma$ is a tree. Choose an orientation $\G0$ of $\Gamma$ 
such that each vertex of $\Gamma$ is either a source or a sink. Equivalently, 
$\G0$ has no oriented paths of length $2.$ The graph $\Gamma$ has exactly 
two such orientations. Notice that $\Bred(\G0)$ 
is isomorphic to the path algebra
of $\G0,$ since the later does not contain any paths of length greater than 
$1.$ 

Below, when we talk about correspondences between indecomposable 
representations, we actually mean correspondences between isomorphism classes 
of indecomposable representations. For brevity, ``isomorphism classes'' 
will be omitted everywhere. 

\begin{prop} There is a natural two-to-one 
correspondence between indecomposable representations of $A(\Gamma)$ and 
indecomposable representations of $B(\G0).$ 
\end{prop} 

\emph{Note: } The referee pointed out that this is a known result, 
proved in Tachikawa \cite{Ta}. We retained the proof for completeness.  

\emph{Proof: } 
Let $M$ be a representation of $A(\Gamma).$ If for 
some vertex $a$ of $\Gamma$ the module $M$ is not annihilated by $(a|b|a),$ 
then it is easy to see that $M$ contains a projective module as a direct 
summand. Therefore, we assume from now on that $M$ is annihilated by 
$(a|b|a)$ for each vertex $a.$ The module $M$ decomposes as a direct sum 
of vector spaces $(a)M,$ over all $a.$ 
Denote by $v_+$ the set of vertices of $\Gamma$ which are source vertices 
of $\Gamma^0$ and by $v_-$ the set of sinks of $\Gamma^0.$ The disjoint 
union of $v_+$ and $v_-$ is the set of vertices of $\Gamma.$ 

For each edge $(a,b)$ of $\Gamma$ there are maps $(a)M \to (b)M$ and 
$(b)M\to (a)M,$ given by left multiplication by $(b|a)$
and $(a|b),$ respectively. For each $a$ denote by $Ker(a)$ the subspace of 
$(a)M$ which is the intersection of kernels of all maps $(a)M\to (b)M,$ 
as $b$ varies over all vertices adjoint to $a.$ For each $a$ choose 
an arbitrary complement $Comp(a)$ to $Ker(a)$ in $(a)M.$ Let 
\begin{eqnarray*} 
 M_0 &= & \oplusop{a\in v_+}Comp(a)\oplus (\oplusop{a\in v_-}Ker(a)) \\
 M_1 &= & \oplusop{a\in v_-}Comp(a)\oplus (\oplusop{a\in v_+}Ker(a)) 
\end{eqnarray*}
Since $M$ is annihilated by paths $(a|b|a),$ for all pairs $(a,b),$ 
vector subspaces $M_0$ and $M_1$ are actually submodules of $M.$
Submodule $M_0$ has the following structure: we start with a $B(\Gamma^0)$ 
module and extend the action to the whole $A(\Gamma)$ by declaring that
the standard complement $B(\Gamma^0)^{\ast}$ to $B(\Gamma^0)$ in $A(\Gamma)$
acts by $0$. In particular, $M_0$ is  indecomposable if and only if 
it is indecomposable as a $B(\Gamma^0)$-module. Similarly, $M_1$ is 
indecomposable if and only if it is indecomposable as a $B(\Gamma^1)$-module. 

We thus get a map from the disjoint union of indecomposable representations
of $B(\Gamma^0)$ and $B(\Gamma^1)$ to the set of indecomposable 
representations of $A(\Gamma).$ This map restricts to a bijection from 
the set of non-irreducible indecomposable representations of 
 $B(\Gamma^0)$ and $B(\Gamma^1)$ to the set of indecomposable 
$A(\Gamma)$ representations which are neither irreducible nor projective. 
On irreducible representations this map is two-to-one. Since there is a 
one-to-one correspondence between irreducible and projective representations
of $A(\Gamma),$ we can modify this map to be a bijection between 
the disjoint union of indecomposable representations
of $B(\Gamma^0)$ and $B(\Gamma^1)$ to the set of indecomposable 
representations of $A(\Gamma).$ There is a canonical modification which 
preserves the symmetry between $\Gamma^0$ and $\Gamma^1.$ Namely, if 
$a\in v_+,$ send the simple $B(\Gamma^0)$ module $\C(a)$ to the projective 
$A(\Gamma)$-module $P_a$ and the simple $B(\Gamma^1)$-module $\C(a)$ to 
the simple quotient of $P_a.$ For $a\in v_-,$ 
 send the simple $B(\Gamma^1)$ module $\C(a)$ to the projective 
$A(\Gamma)$-module $P_a$ and the simple $B(\Gamma^0)$-module $\C(a)$ to 
the simple quotient of $P_a.$

To finish the proof, note that there is a bijection 
between indecomposable representations of $B(\Gamma^0)$ and 
$B(\Gamma^1)$ given by passing to the dual vector space. $\square$

\vspace{0.1in}

For a finite Dynkin graph $\Gamma$ indecomposable representations of 
$B(\G0)$ are in a one-to-one correspondence with positive roots 
of the root system associated to $\Gamma.$ We get a 

\begin{corollary} 
If $\Gamma$ is  a finite Dynkin graph then  indecomposable representations of 
$A(\Gamma)$ are in a one-to-one correspondence with roots of $\Gamma.$ 
\end{corollary} 

\section{Zigzag algebras for affine Dynkin diagrams and the McKay 
correspondence} 
\label{McKay} 

\subsection{Zigzag algebras, Koszulity and quantum Cartan matrices} 
\label{q-Cartan-mat} 

If $\Gamma$ has more than one vertex, $A(\Gamma)$ is 
generated by elements of degree $0$ and $1.$ If $\Gamma$ has more than 
two vertices, $A(\Gamma)$ is a quadratic algebra. 
Roberto Mart\'{\i}nez-Villa [MV] 
found a surprising characterization of Dynkin diagrams in 
terms of zigzag algebras $A(\Gamma)$: 

\begin{prop} Algebra $A(\Gamma)$ is Koszul if and only if $\Gamma$ is 
not a finite Dynkin graph. 
\end{prop} 

The quadratic dual $A^!(\Gamma)$ of $A(\Gamma)$ is isomorphic to the 
quotient algebra of the path algebra of the double graph $D\Gamma$ 
by relations $\sum_b (a|b|a)=0$ where we sum over all vertices $b$ of
$\Gamma$ adjacent to $a.$ If $\Gamma$ is bipartite, $A^!(\Gamma)$ 
is isomorphic to the preprojective algebra of the source-sink oriented
graph $\Gamma$ (see Reiten [R] for an introduction to preprojective algebras). 

Note that the algebra  $A^!(\Gamma)$ is finite-dimensional 
if and only if $\Gamma$ is a finite Dynkin diagram, while 
$A(\Gamma)$ is finite-dimensional for any graph $\Gamma.$

The Cartan matrix of  a finite-dimensional $\C$-algebra $B$ has rows and 
columns enumerated by isomorphism classes of indecomposable projective 
left $B$-modules, and its entries are dimensions of homomorphism spaces
between these projectives. Algebra $A(\Gamma)$ is graded, and its 
Cartan matrix has coefficients which are polynomials in $q,$ 
\begin{equation} 
c_{ab} = \sum_{i\ge 0} q^i \mbox{dim}(\mbox{Hom}_{A(\Gamma)\mbox{-Mod}}
(P_a\{ i\}, P_b)),  
\end{equation} 
where $a$ and $b$ are vertices of $\Gamma.$ Clearly, $c_{a,a}= 1+ q^2,$ 
$c_{a,b}= q$ if $a$ and $b$ are connected by an edge and $c_{a,b}=0$ 
otherwise. 

On the other hand, 
in the theory of Lie algebras the expression 
"Cartan matrix" is 
used to denote a matrix naturally associated to a graph $\Gamma.$ 
This matrix has $2$ as each diagonal entry, $-1$ on the intersection 
of the column $a$ and row $b$ if $a$ and $b$ are joined by an edge 
and $0$ otherwise. 
Now, if we set $q=-1,$ the Cartan matrix of the algebra 
$A(\Gamma)$ is equal to the Cartan matrix of the graph $\Gamma.$ 

The determinant of the Cartan matrix of $\Gamma$ is $0$ if 
$\Gamma$ is an affine Dynkin diagram. However, the 
determinant of the Cartan matrix of $A(\Gamma)$ is a nonzero polynomial 
in $q$ for any graph $\Gamma.$ We will call the Cartan matrix of 
$A(\Gamma)$ {\it the quantum Cartan matrix} of $\Gamma$ and 
denote it by $C(\Gamma).$ This matrix  is always 
invertible in the field $\Q(q)$ of rational functions in $q.$

Let $P^!_a = A^!(\Gamma) (a)$ be the indecomposable projective left 
$A^!(\Gamma)$ module associated to the minimal idempotent $(a)$ of 
$A^!(\Gamma),$ for a vertex $a$ of $\Gamma.$ 
Let $C^!(\Gamma)$ be the Cartan matrix of $A^!(\Gamma).$ Its $(a,b)$-entry 
is 
\begin{equation} 
\label{jab}
c^!_{a,b}= \sum_{i\ge 0} q^i \mbox{dim}(\mbox{Hom}(P^!_a\{ i\}, P^!_b)) 
\end{equation} 
As we have already mentioned, 
if $\Gamma$ is not a finite Dynkin diagram, the zigzag algebra $A(\Gamma)$ 
is Koszul. This and the acyclicity of the Koszul complex implies

\begin{prop}
\label{add-minus} If 
$\Gamma$ is not a finite Dynkin diagram, the matrix $C^!(\Gamma),$ with 
$q$ changed to $-q$ everywhere, is the inverse matrix of $C(\Gamma):$
\begin{equation*}
 C_q(\Gamma) C^!_{-q}(\Gamma) = \mathrm{I}.
\end{equation*}
\end{prop} 

\begin{corollary} If $\Gamma$ if not a finite Dynkin diagram, the 
entries of the inverse quantum Cartan matrix $C(\Gamma)^{-1}$ 
are power series in $-q$ with nonnegative coefficients. 
\end{corollary}

Quantum Cartan matrices for affine Dynkin diagrams $\Gamma$ appear 
in the paper of Lusztig and Tits [LT] which contains a simple algorithm 
for computing the inverse of a Cartan matrix. Since the Cartan matrix 
of an affine Dynkin diagram is not invertible, Lusztig and Tits 
change all diagonal entries from $2$ to $T+ T^{-1}.$ If we set $T=-q$ 
and multiply each entry by $-q,$ we get our quantum Cartan matrix 
$C(\Gamma).$ 

\subsection{A digression: Koszul duality for cross-products}

Let $V$ be a finite-dimensional complex vector space and $G$ a finite group 
acting on $V.$ Let $SV= \oplusop{i\ge 0} S^iV$ be the polynomial algebra 
of $V$ and $S(V,G)= SV\o \C[G]$ be the cross-product algebra. The 
multiplication in $S(V,G)$ is given by 
\begin{equation}
\label{product}
(a\o g)(b\o h) = a g(b) \o gh, \hspace{0.1in} \mbox{ for } 
\hspace{0.1in} a,b\in SV \mbox{ and } g,h\in G. 
\end{equation} 
$S(V,G)$ is a $\Z$-graded algebra, $S(V,G)= \oplusop{i\ge 0}
S^i(V,G)$ where $S^i(V,G)= S^iV\o \C[G].$
Note that $S^0 (V,G)= \C[G]$ is semisimple and that  $SV$ and $\C[G]$ are 
subalgebras of $S(V,G).$ 

The left $S(V,G)$-module $\C[G]= S(V,G)/ S^{>0}(V,G)$ has a projective 
resolution by left $S(V,G)$-modules 
\begin{equation}
\label{resolution}  
\dots \lra SV\o \Lambda^2V \o \C[G]\lra SV\o V \o \C[G] \lra \C[G]\lra 0
\end{equation} 
where the differential 
$\partial: SV \o \Lambda^i V \o \C[G]\lra SV \o \Lambda^{i-1} V \o \C[G]$
is
\begin{equation} 
\label{differential} 
\partial(x\o y_1 \dots y_i \o g) = \sum_{j=1}^i (-1)^j xy_j \o y_1\dots 
y_{j-1}y_{j+1} \dots y_i \o g
\end{equation} 
and the left $S(V,G)$-module structure of $SV\o \Lambda^i V\o \C[G]$ 
is given by 
\begin{equation} 
(x_1\o g)(x_2 \o y \o h) = x_1 g(x_2) \o g(y) \o gh, 
\hspace{0.15in} x_1,x_2\in SV, \hspace{0.05in} 
y\in \Lambda^i V, \hspace{0.05in}g,h\in G. 
\end{equation} 

It is easy to see that $SV\o \Lambda^i V \o \C[G]$ is a free 
$S(V,G)$-module of rank equal to the dimension of $\Lambda^i V.$ 
In particular, (\ref{resolution}) is a graded projective resolution and 
its $i$-th term is generated by the subspace $1\o \Lambda^i V\o 1$ of 
elements of degree $i.$ Therefore, we obtain 

\begin{prop} Graded algebra $S(V,G)$ is quadratic and Koszul. 
\end{prop} 

Define $\Lambda(V,G)$ as the cross-product algebra $\Lambda V \o \C[G]$ 
where $\Lambda V$ is the exterior algebra of $V,$ and the 
product in $\Lambda(V,G)$ is also given by the equation (\ref{product}). 
Algebra $\Lambda (V,G)$ has a natural grading coming from the grading 
of the exterior algebra $\Lambda V.$ 

Consider the algebra $\Lambda(V^{\ast},G)$ where $V^{\ast}$ is the 
dual representation of $G.$ The differential (\ref{differential}) 
commutes with  the natural right action of $\Lambda(V^{\ast},G)$ on 
$SV\o \Lambda V\o \C[G]$ (the convolution action of $\Lambda(V^{\ast})$ 
on $\Lambda V$ together with the right multiplication in the group 
algebra), and we can view (\ref{resolution}) 
as the Koszul complex for the pair of Koszul dual algebras $S(V,G)$ and 
$\Lambda(V^{\ast},G).$ 

\begin{corollary} \label{duallity} The graded 
algebra $\Lambda(V^{\ast},G)$ is quadratic and 
Koszul. Its Koszul dual is isomorphic to $S(V,G).$ 
\end{corollary} 

\emph{Remark} Strictly speaking, since the group ring of $G$ is in 
general not 
commutative, we should distinguish between left and right 
quadratic duals (see [BGS], Section 2.8). However, in our 
case, left and right quadratic duals are isomorphic, and that distinction 
is not necessary. 

\emph{Remark} Suitably formulated, the duality of Corollary~\ref{duallity} 
holds for any reductive group $G$ 
and a finite-dimensional representation $V$ of $G.$ 

\begin{prop} \label{assumptions}
Algebra $\Lambda(V,G)$ is Frobenius. If $\mathrm{dim}V$ is
odd and $G\subset SL(V)$ then $\Lambda(V,G)$ is symmetric. If 
$\mathrm{dim}V$ is even, $G\subset SL(V)$ and $G$ contains a central 
element $h$ acting as $-\mathrm{Id}$ on $V$ then $\Lambda(V,G)$ is symmetric. 
\end{prop} 

\emph{Proof:} Let $n=\mathrm{dim}V$ and $x\in \Lambda^n V,x\not= 0.$ 
The trace map $\mathrm{tr}$ which is $0$ on $\Lambda^i V\otimes \C[G]$ for 
$i<n$ and $\mathrm{tr}(x\otimes g)= \delta_{g,1}$ makes $\Lambda(V,G)$ 
Frobenius. This trace is symmetric when $(V,G)$ satisfies the first 
condition of the proposition. If $(V,G)$ satisfies the second condition, 
set instead $\mathrm{tr}(x \otimes g) = \delta_{g,h},$ for $g\in G.$   
$\square$ 

\subsection{Zigzag algebras and resolutions of simple singularities} 
\label{simplesing} 

Let $G$ be a finite subgroup of $SU(2)$ and $X=\widetilde{\C^2/G}$ 
the minimal resolution of the quotient $\C^2/G.$ The singular fiber is
a union of projective lines, to which we associate a graph with one vertex 
for each projective line and an edge for each pair of intersecting projective 
lines. We denote this graph by $\Gamma(G)$, or simply by $\Gamma.$ This 
graph is a finite Dynkin diagram, and the construction above gives 
a well-known bijection between finite subgroups of $SU(2)$ and simple 
simply-laced Lie algebras. 

To catch a glimpse of zigzag algebras, form the direct sum $\mc{O}'$ of 
the structure sheaves of projective lines in $X.$ 

\begin{prop}\label{appetizer} 
Graded algebras $\mathrm{Ext}^{\ast}_{\mathrm{Coh}(X)}
(\mc{O}', \mc{O}')$ and $A(\Gamma)$ are naturally isomorphic. 
\end{prop} 

The Ext groups are computed in the category of coherent sheaves on $X.$ 
This proposition also appears in [ST]. $\square$ 

We will say that $G$ is \emph{binary} if $G$ has even order, equivalently, 
if $G$ contains $-I,$ the only order $2$ element of $SU(2).$ 
A non-binary subgroup is necessarily an odd order cyclic group. 

The action of   $G$ on $\C^2$ naturally extends to an action of $G$ on 
the exterior algebra on two generators. Form the cross-product algebra 
$\Lambda(\C^2, G).$ 
This is a finite-dimensional algebra, Morita 
equivalent to an algebra described by a finite quiver with relations. 
Vertices of this quiver are in a bijection with irreducible representations 
of $G,$ since the zero degree component of $\Lambda(\C^2,G)$ 
is isomorphic to the group algebra of $G.$ 
Denote by $V_a$ the irreducible representation of $G$ 
associated with the vertex $a.$ 
Recall that, as observed by McKay [McK], 
if to a $G\subset SU(2)$ we associate a graph with 
vertices -- irreducible representations of $G$ and with the number of edges 
connecting $a$ and $b$ equal to the multiplicity of $V_b$ in the 
tensor product $\C^2\o V_a,$ we obtain an affine Dynkin diagram. We denote 
this diagram by $\Gamma^{\mathrm{aff}}.$ 

In our case, the oriented graph underlying the quiver algebra of 
$\Lambda(\C^2,G)$ is the oriented double of $\Gamma^{\mathrm{aff}}.$ 
Since

(a)  $\Lambda(\C^2,G)$ has the  top degree component in degree two, 

(b) if $G$ is binary, $\Lambda(\C^2, G)$ has a symmetric nondegenerate 
 graded trace (Proposition~\ref{assumptions}), 

we easily deduce 

\begin{prop} \label{m-equivalent} If $G$ is binary, 
$\Lambda(\C^2, G)$ is Morita equivalent to the 
zigzag algebra of the affine Dynkin diagram $\Gamma^{\mathrm{aff}}.$ 
\end{prop} 

$G$ is non-binary iff it is cyclic of odd order. In general, 
if $G$ is cyclic of order $n$ then  
$\Lambda(\C^2,G)$ is isomorphic to the quiver algebra of the quiver 
with vertices $1,2, \dots ,n,$ edges $(i|i\pm 1)$ modulo $n,$ and 
relations 
\begin{equation}
\label{skew-z}
(i|i+1|i+2)=(i|i-1|i-2)=0, \hspace{0.2in} (i|i-1|i)+(i|i+1|i)= 0.
\end{equation}
This is an example of a skew-zigzag algebra (see Section~\ref{braid-act}). 
It is isomorphic to a zigzag algebra if $n$ is even.

The Koszul dual statement to Proposition~\ref{m-equivalent} can be 
formulated as  

\begin{prop} If $G$ is binary, 
$S(\C^2, G)$ is Morita equivalent to $A^!(\Gamma^{\mathrm{aff}}),$ the 
latter isomorphic to the preprojective algebra of $\Gamma^{\mathrm{aff}}.$  
\end{prop} 

Kapranov and Vasserot [KV] proved that the derived category of coherent 
sheaves on the minimal resolution $X$ is equivalent to the derived 
category of finitely-generated $S(\C^2,G)$-modules. 
 Notice that $X$ comes with 
a canonical action of $\C^{\ast}.$ We conjecture that the 
derived category of $\C^{\ast}$-equivariant sheaves on $X$ is 
equivalent to the derived category of 
 graded finitely-generated $S(\C^2,G)$-modules. 
If true, then, in view of the Koszul duality between 
$S(\C^2,G)$ and $\Lambda(\C^2,G),$ we would get an equivalence of 
categories between the derived category of $\C^{\ast}$-equivariant sheaves 
on $X$ and the derived category of graded $A(\Gamma^{\mathrm{aff}})$ modules
(Proposition~\ref{appetizer} picks up part of this equivalence).  

\vspace{0.1in}

\begin{prop} The affine braid group associated to the affine Dynkin diagram 
$\Gamma^{\mathrm{aff}}$ acts in the derived category of 
$\Lambda(\C^2, G)$-modules, in the derived category of $S(\C^2,G)$-modules, 
and in the derived category of coherent sheaves on 
the minimal resolution of $\C^2/G.$
\end{prop} 

\emph{Proof:} For $\Lambda(\C^2,G)$ and binary $G$ this follows from 
propositions~\ref{m-equivalent} and \ref{br-action}. For non-binary 
$G$ the algebra $\Lambda(\C^2,G)$ is skew-zigzag and there is a 
braid group action in the derived category according to the remarks 
at the end of Section~\ref{braid-act}. For $S(\C^2,G)$ 
we get the braid group action by Koszul duality. For a vertex $a\in 
\Gamma^{\mathrm{aff}}$ the braid $\sigma_a$ acts by taking 
$M\in D^b(S(\C^2,G)\mbox{-mod})$ to the cone of the evaluation 
map of complexes 
\begin{equation*} 
V'_a\otimes \mbox{RHom}(V'_a, M)\longrightarrow M, 
\end{equation*}
where $V'_a$ is the irreducible $G$-module $V_a,$ considered as an 
$S(\C^2,G)$-module with the trivial action of $S^{>0}(\C^2,G).$ 
Kapranov-Vasserot \cite{KV} equivalence of derived categories 
of $S(\C^2,G)$-modules and coherent sheaves on the minimal resolution of 
$\C^2/G$ allows us to transfer the braid group action to the latter category. 
$\square$ 

Assume that $G\subset SU(2)$ is 
\emph{binary} and let $G'\cong G/\{ \pm 1\}$ be the image of $G$ in $SO(3).$ 
$G'$ acts on $\P^1$ and this action induces an action on the 
cotangent bundle $\TP1.$ Consider the category 
$Coh_{G'}(\TP1)$ of $G'$-equivariant coherent sheaves 
on $\TP1$ and its derived category $D_{G'}(\TP1).$ 

$G'$ acts on $\TP1$ with isolated singular points only. 
The quotient variety $\TP1/G'$ has two or three singular
points (two if $G'$ is cyclic), of type $\C^2/\Z_k,$ where 
$\Z_k\subset SU(2)$ is the stabilizer subgroup of the 
corresponding point on $\P^1.$ 

\begin{prop} The minimal resolution of $\TP1/G'$ is isomorphic to 
the minimal resolution of $\C^2/G.$ 
\end{prop} 

This is proved in Lamotke \cite{La}. $\square$ 

Kapranov-Vasserot theorem [KV] implies 

\begin{prop} The categories of $G'$-equivariant coherent sheaves 
on $\TP1$ and coherent sheaves on the minimal resolution of simple 
singularity $\C^2/G$ are derived equivalent.  
\end{prop} 

\vspace{0.1in} 

What are the multiplicities of various irredicible representations of 
$G\subset SU(2)$ in the $n$-th symmetric power $S^n(\C^2)$ of the 
``basic'' representation $\C^2?$ This problem was solved in different 
ways by Gonzalez-Springer and Verdier [G-SV], Kn\"orrer [Kn],  
Kostant [Ks] and Springer [Sp]. We offer an interpretation
via quantum Cartan matrices (Section~\ref{q-Cartan-mat}) and the  
algebra $S(\C^2, G)$ as follows. 

We restrict to the case of binary $G.$ 
As before, let $a$ be a vertex of $\Gamma^{\mathrm{aff}}$ 
and $V_a$ the irreducible representation of $G$ associated to $a.$ 
Let $t$ be the vertex associated to the trivial representation $V_t$ of $G.$ 
There is an isomorphism of vector spaces 
\begin{equation} 
\label{mod-and-bimod}
\mbox{Hom}_{\C[G]}(V_a, S^n(\C^2)) \cong 
\mbox{Hom}_{\C[G]^{\o 2}}(V_a\o V_t, S^n(\C^2)\o \C[G])
\end{equation} 
where $g\o h\in \C[G]^{\o 2}$ takes $v_1\o v_2\in V_a\o V_t$ 
to $g v_1\o v_2$ and $s\o f\in S^n(\C^2)\o \C[G]$ to $gs\o gfh^{-1}.$ 
To $\alpha\in \mbox{L.H.S.}$ the isomorphism associates $\alpha\otimes \beta,$ 
where 
\begin{equation*}
\beta: V_t \lra \C[G], \hspace{0.2in} \beta(v_2)= \sum_{g\in G} g.
\end{equation*}

Notice that $S^n(\C^2)\o \C[G]$ is just the $n$-th degree component 
of the algebra $S(\C^2,G),$ Morita equivalent 
to $A^!(\Gamma^{\mathrm{aff}}).$ We can rewrite the R.H.S. of 
(\ref{mod-and-bimod}) in terms of $A^!(\Gamma^{\mathrm{aff}}).$ 
Namely, $\C[G],$ the degree $0$ component of $S(\C^2,G),$ 
becomes the degree $0$ component of $A^!(\Gamma^{\mathrm{aff}}),$ 
isomorphic to the direct sum of $\C$'s, one for each vertex of the 
affine diagram. Denote this algebra by $A$ and its simple modules 
corresponding to $a$ and $t$ by $L_a,L_t.$ Then 
\begin{equation*}
\mbox{Hom}_{\C[G]^{\o 2}}(V_a\o V_t, S(\C^2)\o \C[G])\cong
\mbox{Hom}_{A\otimes A^{op}}(L_a\otimes L_t,A^!(\Gamma^{\mathrm{aff}})).
\end{equation*} 

Therefore, the right hand side of (\ref{mod-and-bimod}) is isomorphic to 
the $n$-th degree component of the vector space 
$(a) A^!(\Gamma^{\mathrm{aff}}) (t),$ where 
we multiplied $A^!(\Gamma^{\mathrm{aff}})$ on the left, resp. right, 
by minimal idempotents associated to $a,$ resp. $t.$ The dimension 
of this vector space is equal to the coefficient of $q^n$ in the 
power series $c^!_{a,t}$ (see formula (\ref{jab})). Applying
Proposition~\ref{add-minus} we obtain

\begin{prop} If $G$ is binary, the multiplicity of the irreducible 
representation $V_a$ of $G$ in $S^n(\C^2)$ is equal to $(-1)^n$ times 
the coefficient at $q^n$ of the 
$(a,t)$-entry of the inverse quantum Cartan matrix of 
the affine Dynkin diagram associated to $G$ via the McKay correspondence.  
\end{prop} 

\vspace{0.1in} 

\emph{Remark} The restriction on $G$ is unnecessary. 
Essentially the same proof, with the skew-zigzag algebra (\ref{skew-z}) 
substituted for $A^!(\Gamma^{\mathrm{aff}}),$ works when $G$ is 
non-binary. We leave the details to the reader. Numerically, however, the 
case of cyclic $G$ is boring, since every irreducible representation 
is one-dimensional.  

\vspace{0.1in}

We conclude this advertisement of cross-products $S(\C^2,G)$ and 
$\Lambda(\C^2,G)$ as vital ingredients in the McKay correspondence by 
referring the reader to Auslander [Au] and Reiten [R] 
for a relation between $S(\C^2, G)$ and AR quivers of ADE singularities.

\section{Zigzag algebras in representation theory and geometry} 
\label{rep-theory} 

\subsection{Zigzag algebras and finite groups} 
\label{fin-groups} 

Until now we worked over complex numbers and $A(\Gamma)$ was defined over 
$\C$. In fact, $A(\Gamma)$ is defined over the ring of integers, and we 
denote this $\Z$-algebra by $A_{\Z}(\Gamma).$ For a commutative ring $k$ 
denote by $A_k(\Gamma)$ the $k$-algebra $A_{\Z}(\Gamma) \o_{\Z} k.$ 
The following examples show that $A_k(\Gamma)$ and Morita 
equivalent algebras often appear as direct summands of group algebras 
$k[G]$ in finite characteristic. 
Let us denote by $\Gamma_n$ the chain with $n$ vertices: 

\begin{center} \epsfig{file=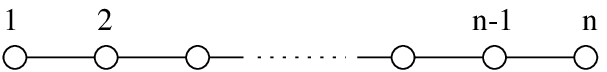} 
\end{center} 

\vspace{0.1in} 

{\bf 1.} Let $p$ be a prime and $k$ a field of characteristic $p.$ 

\begin{prop} The group algebra $k[\SS_p]$ of the symmetric group is the 
direct sum of a semisimple algebra and an algebra Morita equivalent 
to $A_k(\Gamma_{p-1}).$ 
\end{prop} 

This is an exercise in the modular representation theory 
of symmetric groups. $\square$ 

A similar statement holds for Hecke algebras (see Yamane [Y]): 

\begin{prop} The Hecke algebra $H_q(\SS_n)$ when $q$ is an $n$-th primitive 
root of unity is isomorphic to the direct sum of a semisimple algebra 
and an algebra Morita equivalent to $A_{\C}(\Gamma_{n-1}).$ 
\end{prop}

\vspace{0.1in} 

{\bf 2.} 
The group algebra of the finite group $SL(2,p)$ over an algebraically 
closed field $k$ of characteristic $p>2$ decomposes as a direct sum of 
$3$ blocks, one of which is simple and the other two have Cartan 
matrices (see Alperin [Al], \S 17) 
\[ \left(  \begin{array}{ccccccc}
  2 & 1  &       &       &      &      &  \\
  1 & 2  &  1    &       &      &      &  \\
    & 1  &  2    & \cdot &      &      &  \\
    &    & \cdot & \cdot &   1  &      &  \\ 
    &    &       &   1   &   2  &  1   &  \\ 
    &    &       &       &   1  &  2   & 1\\ 
    &    &       &       &      &  1   & 3
\end{array}  \right)   \] 
In each of these two blocks, if we throw away the indecomposable projective 
module $P$ with $\mbox{dim}(\mbox{Hom}(P,P))=3,$ the endomorphism algebra 
$\mbox{End}(P_1\oplus \dots \oplus P_{\frac{p-3}{2}})$ 
of the direct sum of the remaining 
indecomposable projectives is isomorphic to the algebra 
$A_k(\Gamma_{\frac{p-3}{2}}).$ 
The presence of the extra projective $P$ does not create any obstacles 
for defining a braid group action, and  
in the derived category of this block there is a faithful 
action of the braid group on $\frac{p-1}{2}$ strands. 

\vspace{0.1in} 

{\bf 3.}
Let $\mbox{Di}_6$ be the dihedral group of order $12$ and $k$ 
an algebraically 
closed field of characteristic $3.$ Then the group algebra $k[\mbox{Di}_6]$ 
is isomorphic to the direct sum  $A_k(\Gamma_2) \oplus A_k(\Gamma_2).$ 
This result can be easily 
derived from the computation in Curtis and Reiner [CR], \S 91 
 of the Cartan matrix of $\mbox{Di}_6$ in characteristic $3.$ 

\vspace{0.1in} 

{\bf 4.} The following example generalizes $\bf 1,2$ and $\bf 3$. Let 
$B(n,r)$ be the quiver algebra 
\begin{center} \epsfig{file=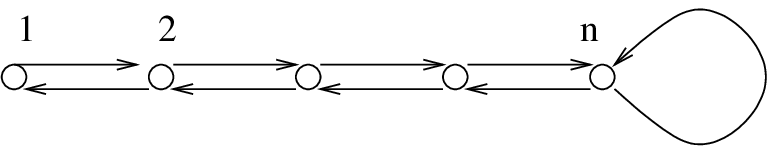} 
\end{center} 
with defining relations 
\begin{eqnarray*} 
& & (i|i+1|i+2)=(i|i-1|i-2)=(n-1|n|n)=(n|n|n-1)=(n|n)^{r+1}=0,  \\
& & (i|i-1|i)=(i|i+1|i), \hspace{0.1in} (n|n-1|n)=(n|n)^r.
\end{eqnarray*} 

Algebra $B(n,r)$ is an example of a Brauer tree algebra with $n+1$ vertices 
and exceptional multiplicity $r.$ Alperin~\cite{Al} calls it  
the \emph{open polygon} Brauer tree algebra (vertices in Brauer trees 
do not correspond to vertices in the quiver algebra, hence the discrepancy 
between $n$ and $n+1$). If $r>1,$ we will call the rightmost vertex 
exceptional. 

Brauer tree algebras and Morita equivalent algebras are isomorphic to 
the so-called 
cyclic defect blocks of group algebras of finite groups over algebraically 
closed fields of finite characteristic (\cite{Al},\cite{Fe},\cite{KZ}). 
This description of cyclic defect blocks was an important achivement of the 
modular representation theory. Cyclic defect blocks are quite widespread. 
For example, the classification of  all blocks of cyclic defect and 
their Brauer 
trees in sporadic simple groups takes up over four hundred pages in the 
monograph \cite{HL}. 

Rickard (\cite{Ri1}, see also \cite{KZ}, chapters 5 and 10) 
showed that any Brauer tree algebra 
with $n+1$ vertices and multiplicity $r$ is derived equivalent to 
$B(n,r).$ 

The subalgebra of the algebra $B(n,r)$ consisting of all paths that neither 
start nor end in the exceptional vertex is isomorphic to the zigzag 
algebra $A(\Gamma_{n-1})$ (or $A(\Gamma_n)$ if $r=1$).

In particular, the results of \cite{KS} imply that there is a faithful 
action of the $n$-stranded braid group in the derived category of 
$B(n,r)$ and of the $(n+1)$-stranded braid group in the derived 
category of $B(n,1).$ The braid group generators act  
by tensoring with the complex of bimodules 
\begin{equation*} 
0 \lra P_i \o \hspace{0.03in} _iP \lra B(n,r) \lra 0, 
\end{equation*}
where $P_i,$ respectively $_iP,$ is the left, resp. right, indecomposable 
projective for the vertex $i.$ 

\vspace{0.1in} 

When this paper was ready for publication we learned 
that this braid group action was independently discovered by A.Zimmermann and  
R.Rouquier in \cite{RZ}. They proved that the braid group action 
is faithful for the $3$-stranded braid group. In addition, the survey paper of 
Rouquier \cite{Ro} offers several startling conjectures about braid 
group actions in derived categories. 

\vspace{0.1in} 

{\bf 5.}
Let $\Gamma$ be the cyclic graph with $3$ vertices and $3$ edges and 
$k$ a field of characteristic $2$ that contains a cubic root of $1.$ 
The algebra 
$A_k(\Gamma)$ is isomorphic to the group algebra (over $k$) of the 
alternating group $\mathbb A_4$ (see Erdmann [E], page 62). 
More examples of blocks Morita equivalent and derived Morita equivalent 
to $A_k(\Gamma)$ can be extracted from Section 12.6 of Feit \cite{Fe}
(and see Rickard \cite{Ri2} for a derived equivalence between the group 
algebra of $\mathbb A_4$ and the principal block 
of $\mathbb A_5$ in characteristic $2$).  

\vspace{0.2in} 

\subsection{Zigzag algebras and the Lie algebra $\mf{sl}_2$} 
\label{and-sl2} 

Let $V$ be the fundamental representation of the Lie algebra $\mf{sl}_2$
over $\C.$ Denote by $\mf{sl}_2(V)$ the $5$-dimensional Lie algebra 
isomorphic as a vector space to $\mf{sl}_2\oplus V$ with the 
Lie bracket 
\begin{equation} 
\label{bracket-plus} 
[(x,a), (y,b)] = ([x,y], xb- ya), \hspace{0.1in} \mbox{ where } 
 \hspace{0.1in} x,y\in \mf{sl}_2, \hspace{0.05in} a,b\in V. 
\end{equation} 
The following observation is due to Loupias [Lp] 

\begin{prop}\label{longchain} 
The category of finite-dimensional $\mf{sl}_2(V)$ 
representations is equivalent  to the category of  
finite-dimensional modules over the quiver algebra 
\begin{equation} 
\label{sl2} 
 \stackrel{0}{\circ} \doublemaprights{30}{30}{}{} \stackrel{1}{\circ}  
  \doublemaprights{30}{30}{}{} \stackrel{2}{\circ}  
 \doublemaprights{30}{30}{}{} \dots 
   \doublemaprights{30}{30}{}{}
\stackrel{i-1}{\circ}
 \doublemaprights{30}{30}{}{} \stackrel{i}{\circ} 
\doublemaprights{30}{30}{}{} \dots 
\end{equation} 
with relations $(0|1|0)=0$ and $(i|i+1|i)= (i|i-1|i).$ 
\end{prop} 

\vspace{0.1in}

Denote by  $\mf{sl}_2^-(V)$ the $(3|2)$-dimensional super Lie algebra 
 isomorphic as a super vector space to $(\mf{sl}_2, V)$ with the 
super Lie bracket 
\begin{equation} 
[(x,a), (y,b)] = ([x,y], xb+ ya), \hspace{0.1in} \mbox{ where } 
 \hspace{0.1in} x,y\in \mf{sl}_2, \hspace{0.05in} a,b\in V. 
\end{equation}

\begin{prop}\label{verylongchain} 
The category of finite-dimensional $\mf{sl}_2^-(V)$ 
representations is equivalent  to the category of  
finite-dimensional modules over the quiver algebra 
\begin{equation} 
\label{sl2-} 
 \stackrel{0}{\circ} \doublemaprights{30}{30}{}{} \stackrel{1}{\circ}  
  \doublemaprights{30}{30}{}{} \stackrel{2}{\circ}  
 \doublemaprights{30}{30}{}{} \dots 
   \doublemaprights{30}{30}{}{}
\stackrel{i-1}{\circ}
 \doublemaprights{30}{30}{}{} \stackrel{i}{\circ} 
\doublemaprights{30}{30}{}{} \dots 
\end{equation} 
with relations $(i|i+1|i+2)=0, (i|i-1|i-2)=0$ and  $(i|i+1|i)= (i|i-1|i).$ 
\end{prop} 

Note that algebras described in propositions \ref{longchain} 
and \ref{verylongchain}  
are quadratic and in fact quadratic dual to each other and Koszul. 
The second algebra is isomorphic to the zigzag algebra of  
the infinite in one direction chain: 

\begin{center} \epsfig{file=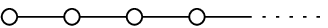} 
\end{center} 

\vspace{0.1in} 

This duality between representations of $\mf{sl}_2(V)$ and 
$\mf{sl}_2^-(V)$ can be generalized to arbitrary pairs $(\mf{g}, V)$ 
where $\mf{g}$ is a semisimple Lie algebra and $V$ a finite-dimensional 
representation of $\mf{g}.$ 
One can form the Lie algebra $\mf{g}(V)= \mf{g}\oplus V$ 
with the bracket (\ref{bracket-plus}) and the "dual"   
Lie superalgebra $\mf{g}^-(V^{\ast}).$ The categories of finite-dimensional 
representations of $\mf{g}(V)$ and $\mf{g}^-(V^{\ast})$ 
are described by Koszul dual algebras.

\subsection{Other examples} 
\label{and-else} 

\vspace{0.1in} 

{\bf 1.} 
Let $T$ be an elliptic curve and $p$ a point of $T.$ Let 
$\mc{L}$ be the direct sum of  
the structure sheaf $\mc{O}_T$ and the skyscraper sheaf $\C_p.$  
 The ext algebra 
$\mbox{Ext}_{\mbox{Coh}(T)}(\mc{L}, \mc{L})$ 
is isomorphic to $A(\Gamma_2).$ Seidel and Thomas~\cite{ST} and 
Thomas~\cite{T} list many other 
appearances of algebras $A(\Gamma)$ as ext algebras of sheaves on 
Calabi-Yau varieties.   

\vspace{0.1in} 

{\bf 2.} 
Algebra $A(\Gamma_n)$ is isomorphic to the algebra of Floer homology 
groups $\oplusop{1\le i,j\le n}HF(L_i, L_j)$ where $L_1, \dots , L_n$ 
is a chain of Lagrangian spheres in a suitable symplectic manifold 
(see [KS]). 
 
\vspace{0.1in} 

{\bf 3.} Let $PS(\mathbb{P}^n)$ be the category of perverse sheaves on 
$\mathbb{P}^n,$ smooth along the stratification of $\mathbb{P}^n$ by 
the increasing chain $\mathbb{P}^0\subset \mathbb{P}^1\subset 
\dots \subset \mathbb{P}^{n}.$ 
Let $Q_i$ be the indecomposable projective perverse 
sheaf associated to the $i$-dimensional strata and 
$Q= \oplusop{0\le i\le n-1}Q_i.$ 
The zigzag algebra $A(\Gamma_n)$ is 
isomorphic to the endomorphism algebra $\mbox{End}_{PS(\mathbb{P}^n)}(Q)$ 
(see [KS] for details).

\end{document}